\documentclass[a4paper, 11pt]{article}

\usepackage{amsfonts, amsmath, amsthm, amssymb, stmaryrd, mathtools}
\usepackage[T1]{fontenc}
\usepackage[utf8]{inputenc}
\usepackage[backend=biber,bibstyle=authoryear,citestyle=authoryear-icomp,isbn=false,url=false,doi=false]{biblatex}
\usepackage{latexsym}
\usepackage{ellipsis}
\usepackage{graphicx}
\usepackage{enumitem}
\usepackage{aliascnt}
\usepackage[hidelinks]{hyperref}
\usepackage{cleveref}

\AtEveryBibitem{\clearlist{language} \clearfield{note} \clearfield{edition}}
\usepackage{biblatex}
\renewbibmacro{in:}{%
  \ifentrytype{article}{}{\printtext{\bibstring{in}\intitlepunct}}}

\newtheoremstyle{break}%
{}{}%
{\itshape}{}%
{\bfseries}{}
{\newline}{}

\theoremstyle{plain}
\newtheorem{theorem}{Theorem}
\newtheorem*{theorem*}{Theorem}
\newaliascnt{lemma}{theorem}

\newtheorem*{lemma*}{Lemma}
\aliascntresetthe{lemma}
\newaliascnt{propn}{theorem}
\theoremstyle{plain}

\newtheorem*{proposition*}{Proposition}
\aliascntresetthe{propn}
\newaliascnt{defn}{theorem}
\theoremstyle{definition}

\newtheorem*{definition*}{Definition}
\theoremstyle{break}

\crefname{equation}{}{}
\creflabelformat{equation}{(#2#1#3)}

\renewcommand{\complement}[1]{#1^\mathsf{c}}

\newcommand{\restr}[2]{#1{\restriction}_{#2}}

\begin{document}

\title{A direct proof of Tychonoff's theorem}
\author{Oliver Tatton-Brown}
\date{}
\maketitle

\begin{abstract}
Proofs of Tychonoff's theorem often seem to require a bit of magic. Machinery such as ultrafilters, nets or other devices are employed to give proofs that can be very neat, but which are not the kind of thing that one would naturally think of when presented with the problem (given a background in standard open set topology). Here we present a direct, transparent and pretty simple proof of Tychonoff's theorem, straight from the open cover definition of compactness.
\end{abstract}

Standard proof of Tychonoff's theorem generally seem to require a bit of magic. They employ machinery like ultrafilters, nets or maximal families with the finite intersection property, which are not the kind of thing one would naturally think of when presented with the problem, given a background in standard open set topology \parencites[pp.234--235]{munkres_topology_2000}[pp.138--139]{willard_general_1970}[p.120]{engelking_general_1989}. Here we give a direct, transparent and pretty simple proof of Tychonoff's theorem, straight from the open cover definition of compactness.

First, some basic notation. We write $\text{Dom}(u)$ for the domain of $u$. When thinking of a function $u\in \prod_{b\in B} Y_b$ as an element of the product we write the value of $u$ on $b$ as $u_b$. Restriction of a function to a subdomain is denoted by $\restr{u}{B'}$ as usual.

\begin{theorem}[Tychonoff]
Let $(X_a\mid a\in A)$ be a family of compact topological spaces. Then $X=\prod_{a\in A}X_a$ is compact.
\end{theorem}
\begin{proof}
Let $(U_j\mid j\in J)$ be a family of basic open subsets of $X$ such that if $J'\subseteq J$ is finite, then $\bigcup_{j\in J'}U_j\neq X$. We will show there exists an element $x$ of $X$ such that $x\notin \bigcup_{j\in J}U_j$. That will prove the theorem (since then every family of open subsets of $X$ which does cover $X$ must have a finite subcover of basic open sets).

The proof is by Zorn's lemma. We let $P$ be the set of functions $u\in\prod_{a\in A'}X_a$ for some $A'\subseteq A$ such that for every finite $J'\subseteq J$ there is some $y\in X$ such that $y\in\complement{(\bigcup_{j\in J'}U_j)}$ and $\restr{y}{A'}=u$. $P$ is non empty since it contains the empty function $\varnothing\in\prod_{a\in\varnothing}X_a$, by the initial assumption on the $U_j$.

Before continuing, we note that for each $j$ we can write $U_j$ as $\{y\mid \text{if }a\in A_j\text{ then } y_a\in V^j_a\}$ where $A_j$ is a finite subset of $A$ and each $V^j_a$ is an open subset of $X_a$.

We order $P$ by extension of functions. Then $P$ is easily seen to be chain complete. Indeed, let $C$ be a chain in $P$. Let $u$ be the function $\bigcup_{v\in C}v$, with $\text{Dom}(u)=\bigcup_{v \in C}\text{Dom}(v)$. We will show that $u\in P$. Suppose we have a finite subset $J'$ of $J$. We have $U_j=\{y\mid \text{if }a\in A_j\text{ then } y_a\in V^j_a\}$ with $A_j$ finite for each $j$. Thus $(\bigcup_{j\in J'}A_j)\cap \text{Dom}(u)$ is a finite subset of $\bigcup_{v\in C}\text{Dom}(v)$, so is a subset of $\text{Dom}(v)$ for some $v \in C$, since $C$ is a chain. Then since $v\in P$ we have that there is some $y\in \complement{(\bigcup_{j\in J'}U_j)}$ such that $\restr{y}{\text{Dom}(v)}=v$. Define $z$ by $z_a=u_a$ for $a\in \text{Dom}(u)$, and $z_a=y_a$ for $a\notin \text{Dom}(u)$. Now if $a\in \bigcup_{j\in J'}A_j$ and $a\in \text{Dom}(u)$ then $a\in \text{Dom}(v)$ so $z_a=u_a=v_a=y_a$, and if $a \in \bigcup_{j\in J'}A_j$ and $a\notin \text{Dom}(u)$ then by definition $z_a=y_a$. Thus $z$ agrees with $y$ on $\bigcup_{j\in J'}A_j$, so $z\in U_j$ if and only if $y\in U_j$. That holds for all $j\in J'$, so in fact $z\in \complement{(\bigcup_{j\in J'}U_j)}$. But $\restr{z}{\text{Dom}(u)}=u$ and $J'$ was arbitrary, so indeed $u\in P$ as required.

Thus by Zorn's lemma, $P$ has a maximal element. We will call it $x$. Suppose for contradiction that $\text{Dom}(x)\neq A$, and let $a\in A\setminus \text{Dom}(x)$. For $r\in X_a$ let $x^r$ be the function with domain $\text{Dom}(x)\cup \{a\}$ defined by $\restr{x^r}{\text{Dom}(x)}=x$, $x^r_a=r$. Since $x$ is maximal we have that for all $r\in X_a$, $x^r\notin P$. That means that for every $r\in X_a$ there is some finite subset $J_r$ of $J$ such that if $y\in X$ with $\restr{y}{\text{Dom}(x^r)}=x^r$ then $y\in \bigcup_{j\in J_{r}}U_j$. For each $j$ we have $U_j = \{y\mid \text{if }b\in A_j\text{ then } y_b\in V^j_b\}$ with $A_j$ finite. Without loss of generality we may assume that for each $j\in J_r$, $U_j$ is a neighbourhood of some $y$ with $\restr{y}{\text{Dom}(x^r)}=x^r$; thus if $b \in A_j$ and $b\in \text{Dom}(x^r)$ then $x^r_b\in V^j_b$. Let $A_r=\bigcup_{j\in J_r}A_j$, a finite subset of $A$.

Again without loss of generality we may assume $a\in A_j$ for each $j\in J_r$, since we may set $V^j_a$ to be $X_a$ for $j$ with $a\notin A_j$. Now let $W_r$ be $\bigcap_{j\in J_r}V^j_a$. This is an open neighbourhood of $r$ in $X_a$. Suppose that we have some $y$ with $\restr{y}{\text{Dom}(x)}=x$, $y(a)\in W_r$. Then letting $z\in X$ be defined by $z(a)=r$, $\restr{z}{A\setminus\{a\}}=\restr{y}{A\setminus\{a\}}$, we have $z\in U_j$ for some $j\in J_r$; but then since $y$ agrees with $z$ everywhere but at $a$, and $y_a\in W_r\subseteq V^j_a$, we also have $y\in U_j$. Thus $\{y\mid \restr{y}{\text{Dom}(x)}=x,\,y_a\in W_r\}$ is a subset of $\bigcup_{j\in J_r}U_j$.

Every $W_r$ is an open neighbourhood of $r$, so we can find $K\subseteq X_a$ finite with $\bigcup_{r\in K}W_r=X_a$. But then for any $y$ with $\restr{y}{\text{Dom}(x)}=x$, we have $y_a\in W_r$ for some $r\in K$, so by the above we obtain $y\in \bigcup_{j \in J_r}U_j$. Thus letting $J'=\bigcup_{r\in K}J_r$, $J'$ is a finite subset of $J$, such that for every $y$ with $\restr{y}{\text{Dom}(x)}=x$, we have $y\in\bigcup_{j\in J'}U_j$. But that contradicts the fact that $x\in P$.

This means that our assumption, that $\text{Dom}(x)\neq A$, must be false. In other words $\text{Dom}(x)=A$. Now let $j$ be any element of $J$. Since $x\in P$ there is some $y\in X$ with $\restr{y}{\text{Dom}(x)}=\restr{x}{\text{Dom}(x)}$, and $y\notin U_j$. But that just means that $y=x$, so $x\notin U_j$. This holds for all $j$, so we are done.
\end{proof}

\printbibliography

\end{document}